\documentclass[12pt]{a_t}
\usepackage{graphicx}

\year{2010}
\title{К оптимизации линейных систем с управляемыми коэффициентами}%
\thanks{Работа выполнена при частичной финансовой поддержке Российского фонда фундаментальных исследований (грант \mbox{№\,08-08-90026}-Бел\_а).}

\authors{%
О.~В.~Батурина,\\
А.~В.~Булатов, канд.~физ.-мат.~наук,\\
В.~Ф.~Кротов, д-р~техн.~наук\\
(Институт проблем управления им. В.~А.~Трапезникова РАН, Москва)}

\abstract{%
Исследуются характерные свойства линейных дифференциальных уравнений с управляемыми коэффициентами, их динамические инварианты и ограничения области достижимости, существенные для проблем управления, свойства понтрягинских экстремалей данного класса задач, численные методы оптимизации управления. PACS 02.30.Xx, 02.30.Yy, 02.60.Pn
}

\newcommand{\Argmax}{\mathop{\mathrm{Argmax}}}

\newcommand{\Real}{\mathbb{R}}
\newcommand{\argmax}{\mathop{\mathrm{argmax}}}
\newcommand{\transposed}{\intercal}

\begin{document}
\maketitle

\section*{Введение}

Рассматривается  задача оптимального управления динамическими системами, описываемыми линейными дифференциальными уравнениями с управляемыми коэффициентами и квадратичным критерием оптимальности. Эта модель характерна для многих  прикладных задач, играющих существенную роль в теории управления. Это, например, задача синтеза регулятора нелинейных систем, обеспечивающего достаточно малое отклонение траектории возмущенного движения от заданного. Если соответствующие отклонения управления также  малы, то применима линеаризация в окрестности заданного движения и задача сводится к хорошо разработанной линейно-квадратичной задаче  оптимального управления (ЛКЗ). Если же малое отклонение траектории возмущенного движения от заданного достигается при больших отклонениях  управления, то применима лишь частичная линеаризация, и задача сводится к рассматриваемой здесь нелинейной вариационной задаче.

Другой пример: широкий и актуальный класс задач  управления квантовым состоянием атомов и молекул. Динамика таких систем описывается дифференциальным уравнением Шредингера, относящимся к  рассматриваемому классу, а цели управления адекватно формализуются квадратичными критериями оптимальности.

Исследуются характерные свойства линейных дифференциальных уравнений с управляемыми коэффициентами, их динамические инварианты и ограничения областей достижимости, существенные для проблем управления, свойства понтрягинских экстремалей  данного класса задач, численные методы оптимизации управления.

\section{Постановка задачи}\label{MainProblem}

Рассмотрим задачу (далее задача I)
\begin{equation}\label{MainProblem1}
   I(x;u)=(x(T),Lx(T)) \to\min,
\end{equation}
\begin{equation}\label{MainProblem2}
   \dot x=(A+uB)x,
\end{equation}
\begin{equation}\label{MainProblem3}
   x(0)=x_0,
\end{equation}
\begin{equation}\label{MainProblem4}
   |u(t)| \leqslant \nu, \qquad t\in[0,T],
\end{equation}
где $x = (x^1, \ldots, x^n)$ - вещественный фазовый вектор; $L$ - вещественная симметричная матрица; $A$ и $B$ - вещественные матрицы; управление $u(\cdot)$ одномерно; функции $x(\cdot)$ и $u(\cdot)$ определены на отрезке $[0,T]$; числа $T>0$ и $\nu>0$, вектор $x_0$, матрицы $A$, $B$, $L$ заданы. Решение задачи (\ref{MainProblem1})--(\ref{MainProblem4}) понимается в классическом смысле: управление $u(\cdot)$ кусочно непрерывно на $[0,T]$, а траектория $x(\cdot)$ непрерывна и имеет кусочно непрерывную производную.

Используются следующие стандартные обозначения: $(x,y) = \sum_i x^i y^i$  - евклидово скалярное произведение в пространстве $\Real^n$; $\|x\| = [ (x,x) ]^{1/2}$ - евклидова норма вектора; $\|A\|=\left[ \sum_{ij} a_{ij}^2\right]^{1/2}$ - евклидова норма матрицы.

Вместо задачи I часто решают задачу \eqref{MainProblem1Beta}, (\ref{MainProblem2})--(\ref{MainProblem4})  (далее задача I${}_\beta$):
\begin{equation}\label{MainProblem1Beta}
   I_\beta(x;u)=(x(T),Lx(T)) + \beta \int\limits_0^T (u(t))^2 dt \to\min,
\end{equation}
где $\beta>0$ задано. Использование регуляризующего слагаемого в (\ref{MainProblem1Beta}) упрощает исследование уравнений принципа максимума и итерационные алгоритмы оптимизации (исключаются особые режимы управления). Но при этом вносится неясность в физический смысл критерия, обуславливается недоиспользование ресурса управления для главной цели (\ref{MainProblem1}) и ухудшается сходимость итерационных методов. Подробнее об этом см. в \cite{Krotov2009Ait}.

\section{Некоторые свойства уравнений системы и оптимизационной задачи}\label{SectPropsOfSystems}

1. Область достижимости системы (\ref{MainProblem2}), (\ref{MainProblem3}) ограничена и отделена от начала координат. Справедливо следующее утверждение:

\begin{theorem}\label{Th1}
При всех допустимых в задаче I управлениях и при всех{ } $t\in[0,T]$ имеет место двусторонняя оценка
\begin{equation}\label{DomainEstimate}
  \| x_0 \| e^{-\gamma T} \leqslant \|x(t)\| \leqslant   \| x_0 \| e^{\gamma T},
\end{equation}
где
\begin{displaymath}
  \gamma=\|A\| + \nu \|B\|.
\end{displaymath}
\end{theorem}
Доказательство этой теоремы приводится в Приложении.

2. Выделим частный случай системы (\ref{MainProblem2}), (\ref{MainProblem3}): размерность фазового вектора четная ($n=2m$), а $x=(x^1, x^2, \ldots, x^n) = (y^1, y^2, \ldots, y^m, z^1, z^2, \ldots, z^m)$; матрицы $A$ и $B$ имеют следующую блочную структуру:
\begin{equation}\label{BlockStruct}
  A=\begin{pmatrix} 0 & P_A \\ -P_A & 0 \end{pmatrix},
  \qquad
  B=\begin{pmatrix} 0 & P_B \\ -P_B & 0 \end{pmatrix},
\end{equation}
где $P_A$ и $P_B$ - симметричные матрицы размера $m\times m$. Тогда система (\ref{MainProblem2}) представляет собой управляемую гамильтонову систему вида
\begin{equation}\label{DotY}
   \dot y = \frac{\partial {\cal H}}{\partial z} = P z,
\end{equation}
\begin{equation}\label{DotZ}
   \dot z = -\frac{\partial {\cal H}}{\partial y} = - P y,
\end{equation}
где $P=P_A+uP_B$, при этом гамильтониан системы имеет вид
\begin{equation}\label{Hamiltonian}
   {\cal H}(y,z) = \frac{1}{2} (Py,y) + \frac{1}{2} (Pz,z).
\end{equation}
Частный случай задачи I (I${}_\beta$), в котором матрицы $A$ и $B$ имеют блочную структуру (\ref{BlockStruct}), назовем задачей II (соответственно II${}_\beta$).

3. В случае (\ref{BlockStruct}) система (\ref{MainProblem2}) имеет компактное комплексное представление. Рассмотрим комплексное уравнение, представляющее собой конечномерный аналог основного уравнения квантовой механики (уравнения Шредингера в частных производных)
\begin{equation}\label{Schroedinger}
  \dot w = iHw,
\end{equation}
где $i$ - мнимая единица, $H=H^1+iH^2$ - комплексный линейный самосопряженный оператор, $w=y+iz$ - комплексный вектор фазовых переменных. Выделяя в этом уравнении вещественную и мнимую часть получаем систему
\begin{displaymath}
  \dot y = -H^2y-H^1z;
\end{displaymath}
\begin{displaymath}
  \dot z =  H^1y-H^2z.
\end{displaymath}
При $H^2=0$ эта система имеет вид (\ref{DotY}), (\ref{DotZ}).

4. Система (\ref{DotY}), (\ref{DotZ}) имеет динамический инвариант
\begin{equation}\label{Invar}
   f(x) = \|x\|^2 = \|y\|^2 + \|z\|^2 =
   \sum\limits_{k=1}^m (y^k)^2 + \sum\limits_{k=1}^m (z^k)^2 = \mathrm{const}.
\end{equation}
Действительно,
\begin{displaymath}
   \frac{1}{2}\frac{d}{dt} f(x(t)) =
       \sum\limits_{k=1}^m y^k {\dot y}^k + \sum\limits_{k=1}^m z^k {\dot z}^k =
       \sum\limits_{k=1}^m y^k \sum\limits_{s=1}^m p_{ks} z^s -
       \sum\limits_{k=1}^m z^k \sum\limits_{s=1}^m p_{ks} y^s = 0.
\end{displaymath}
Таким образом, все траектории системы (\ref{DotY}), (\ref{DotZ}), независимо от выбора управления, расположены на сфере с центром в начале координат.

Наличие динамического инварианта предоставляет определенную свободу выбора матрицы $L$, фигурирующей в (\ref{MainProblem1}). Так как для любого допустимого управления $\|x(T)\|=\|x_0\|$, то возможна замена исходной матрицы $L$ на матрицу вида $L+\alpha I$, где $\alpha$ - произвольное вещественное число, а $I$ - единичная матрица. В частности, это преобразование позволяет добиться как положительной, так и отрицательной определенности квадратичной формы (\ref{MainProblem1}) в задаче II. Как показано ниже, это обстоятельство существенно упрощает реализацию глобального итерационного метода последовательного улучшения управления. В отличие от задачи II, задача I в общем случае динамического инварианта не имеет.

5. Пусть матрица $L$ является матрицей оператора ортогонального проектирования на подпространство $Q\subset\Real^n$. В этом случае квадратичная форма (\ref{MainProblem1}) выпукла. Пусть $\tilde L$ - оператор ортогонального проектирования на ортогональное дополнение подпространства $Q$. Согласно (\ref{Invar}):
\begin{displaymath}
  (x(T),Lx(T))=\|x(T)\|^2 - (x(T), \tilde Lx(T)) = \|x_0\|^2 - (x(T), \tilde Lx(T)),
\end{displaymath}
и мы имеем пару эквивалентных минимизируемых функционалов: выпуклый (\ref{MainProblem1}) и вогнутый 
\begin{displaymath}
     \tilde I(x;u)=-(x(T),\tilde Lx(T)).
\end{displaymath}
%

\section{Необходимые условия оптимальности (в форме принципа максимума Понтрягина)}

Функция Гамильтона-Понтрягина \cite{Pontryagin} задачи I${}_\beta$ имеет вид
\begin{equation}\label{PontrFunc}
   H(t,\psi, x,u) =
       -\beta u^2 +
      \sum\limits_{i=1}^{n}
          \psi_i \sum\limits_{k=1}^{n} (a_{ik}+ub_{ik})x^k;
\end{equation}
двойственная система
\begin{equation}\label{DualSystem}
    \dot\psi = -\frac{\partial H}{\partial x}  = - (A^\transposed + uB^\transposed)\psi
\end{equation}
($(\cdot)^\transposed$ - операция транспонирования матрицы). В исследуемой задаче двойственная система имеет тот же вид линейной системы с управляемыми коэффициентами, что и исходная. В задаче II сохраняется блочная структура матриц $A$ и $B$:
\begin{displaymath}
  -A^\transposed=\begin{pmatrix} 0 & P_A^\transposed \\ -P_A^\transposed & 0 \end{pmatrix},
  \qquad
  -B^\transposed=\begin{pmatrix} 0 & P_B^\transposed \\ -P_B^\transposed & 0 \end{pmatrix}.
\end{displaymath}
Отсюда следует, что условие инвариантности нормы (п.3 раздела \ref{SectPropsOfSystems}) выполняется не только для исходной системы, но и для двойственной системы.

Условия трансверсальности:
\begin{equation}\label{Transv}
  \psi(T) = - 2Lx(T).
\end{equation}
Управление, отвечающее условию принципа максимума Понтрягина, в задаче I имеет вид:
\begin{equation}\label{ExtremalControl1}
  u(t) = \tilde u(t, \psi(t), x(t)),
\end{equation}
где при $\beta=0$:
\begin{equation}\label{ExtremalControl2}
  \tilde u(t,\psi,x) \in \Argmax\limits_{|u|\leqslant \nu} H(t, \psi, x, u) =
    \left\{
    \begin{array}{cc}
            \nu          &   \text{если ${\cal K}>0$;} \\{}
          \text{любое из } [-\nu,\nu]     &    \text{если ${\cal K}=0$;} \\{}
          -\nu          &    \text{если ${\cal K}<0$,}
    \end{array}
    \right.
\end{equation}
\begin{equation}\label{uterm}
    {\cal K}(\psi, x) = \sum\limits_{i=1}^{n}
          \psi_i \sum\limits_{k=1}^{n} b_{ik}x^k,
\end{equation}
а при $\beta>0$:
\begin{equation}\label{ExtremalControl2Beta}
  \tilde u(t,\psi,x) = \argmax\limits_{|u|\leqslant \nu} H(t, \psi, x, u) =
    \left\{
    \begin{array}{cc}
            -\nu          &   \text{если $\frac{1}{2\beta}{\cal K}(\psi,x) \leqslant -\nu$;} \\{}
          \frac{1}{2\beta}{\cal K}(\psi,x)     &    \text{если $-\nu<\frac{1}{2\beta}{\cal K}(\psi,x)<\nu$;} \\{}
          \nu          &    \text{если $\frac{1}{2\beta}{\cal K}(\psi,x) \geqslant \nu$,}
    \end{array}
    \right.
\end{equation}

Таким образом, для нахождения экстремали задачи I получаем краевую задачу для пары систем (\ref{MainProblem2}), (\ref{DualSystem}) с граничными условиями (\ref{MainProblem3}), (\ref{Transv}), в которой управление связано с фазовыми переменными соотношениями (\ref{ExtremalControl1}), (\ref{ExtremalControl2}). 

Известные методы решения этой краевой задачи сводятся к решению набора задач Коши с фиксированными $\psi_0=\psi(0)$. Значение $\psi_0$ подбирается так, чтобы выполнялись условия трансверсальности (\ref{Transv}). При $\beta>0$ для этих задач Коши выполняются классические условия однозначной разрешимости. При $\beta=0$ возможно существование поверхностей разрыва правой части и особых режимов управления. 

Управление $u_* (\cdot)$ называется {\it особым\/} на отрезке $[t_1, t_2]\subset[0,T]$, если ${\cal K}(\psi(t), x(t)) = 0$ при всех $t\in [t_1, t_2]$.  Здесь $\psi(t)$ и $x(t)$ - траектории двойственной и исходной систем, соответствующие $u_* (\cdot)$. Из (\ref{ExtremalControl1}), (\ref{ExtremalControl2}) следует, что на участке особого режима всегда выполняются условия принципа максимума. Таким образом, оптимальное управление в задаче может содержать участки особого управления. Такие задачи называются вырожденными.

Традиционный способ нахождения особого управления - дифференцирование равенства (\ref{uterm}) по $t$ и подстановка (\ref{MainProblem2}) и (\ref{DualSystem}) в полученное соотношение. Иногда проводить эти операции приходится несколько раз. Для вывода универсального соотношения для особого режима рассмотрим произвольную билинейную форму
\begin{displaymath}
{\cal C} = \sum\limits_{i=1}^{n} \sum\limits_{k=1}^{n} c_{ik} \psi_i x^k.
\end{displaymath}
Вычисляя ее производную по времени с учетом (\ref{MainProblem2}), (\ref{DualSystem}) имеем
\begin{equation}\label{RecEq}
\frac{d{\,\cal C}}{dt} =
 \sum\limits_{k,m}
 \left(
     c_{mk}^{(1)} + u d_{mk}^{(1)}
  \right)
\psi_m x^k,
\end{equation}
где
\begin{displaymath}
C^{(1)} = CA-AC, \quad D^{(1)} = CB-BC.
\end{displaymath}
Для вычисления $d{\cal K}/dt$ положим $C^{(0)}=C=B$, $D^{(0)}=0$. В этом случае $C^{(1)}=BA-AB$, $D^{(1)}=0$. Таким образом, из соотношения $d{\cal K}/dt=0$ особое управление определить не удается: коэффициент при  $u$ обращается в $0$. Проводя дифференцирование билинейной формы $s$ раз (при условии, что $D^{(1)}$ = $D^{(2)}$ = $\cdots$ = $D^{(s-1)}$ = 0), получаем
\begin{displaymath}
\frac{d^s{\,\cal C}}{dt^s} =
 \sum\limits_{k,m}
 \left(
     c_{mk}^{(s)} + u d_{mk}^{(s)}
  \right)
\psi_m x^k,
\end{displaymath}
где
\begin{displaymath}
  C^{(s)} = C^{(s-1)}A - AC^{(s-1)}; \qquad
  D^{(s)} = C^{(s-1)}B - BC^{(s-1)}.
\end{displaymath}
Если при каком-либо $s$ матрица $D^{(s)}$ ненулевая, то из равенства $d^s{\cal K}/dt^s=0$ можно выразить $u$ через $t$ и $x$. Таким образом, получаем соотношение, которому гарантированно удовлетворяет особое управление.

Отметим случай, когда матрицы $A$ и $B$ перестановочны ($AB=BA$). В этом случае $C^{(s)}$ = $D^{(s)}$ = $0$ при всех $s$ и описанная выше процедура не дает возможности найти особое управление. В задаче II такая ситуация всегда имеет место при $m=1$. Оптимальные управления, содержащие один или несколько участков особого режима, встречаются в системах размерности 4 и выше.

Если матрицы $A$ и $B$ перестановочны, то при любом допустимом управлении ${\cal K}(t)$ - константа. В этом случае в силу принципа максимума Понтрягина решением задачи является или управление $u(t)\equiv \nu$ (если для этого управления ${\cal K}>0$) или управление $u(t)\equiv -\nu$ (если для этого управления ${\cal K}<0$) или особое управление. При $n=1$ задача особых управлений не имеет. При $n=2$ задача имеет континуум особых управлений при определенных сочетаниях параметров задачи.

При $t=T$ имеем
\begin{equation}\label{KT}
    {\cal K}(T) =
    \sum\limits_{i=1}^{n} \psi_i(T) \sum\limits_{k=1}^{n} b_{ik}x^k(T) =
    (\psi, Bx) = -2(Lx,Bx)=-2(B^*Lx,x).
\end{equation}
Из (\ref{KT}) следует, что особое управление однозначно характеризуется конечной точкой траектории $x(T)$. Так как равенство $x(T)=0$ в исследуемых задачах невозможно (при $x_0\ne 0$ это противоречило бы единственности решения задачи Коши для системы обыкновенных дифференциальных уравнений), то, в частности, достаточным условием отсутствия особых управлений является положительная или отрицательная определенность квадратичной формы (\ref{KT}).

Из \eqref{ExtremalControl2Beta} следует, что в задаче I${}_\beta$ при $\beta>0$ невозможна ситуация особого режима.

\section{Итерационные методы улучшения управления}

Решение краевой задачи принципа максимума (\ref{MainProblem2}), (\ref{DualSystem}), (\ref{MainProblem3}), (\ref{Transv}), (\ref{ExtremalControl1}) технически возможно лишь для задач небольшой размерности. На практике для решения задач большой размерности (в конкретных задачах квантовой механики размерность фазового вектора имеет порядок $10^4$-$10^6$) используются итерационные методы последовательного улучшения управления: широко известный градиентный метод (см. например \cite{GradMethodGeneral}, \cite[стр. 239]{MDekker1996}), а также глобальный метод \cite[\S {}6.2]{MDekker1996}, имеющий применительно к данной задаче определённые преимущества.

Опишем глобальный метод применительно к общей задаче оптимального управления
\begin{equation}\label{GMEq1}
  I(x;u) = \int\limits_0^T f^0 (t,x,u)dt + F(x(T)) \to\min,
\end{equation}
\begin{equation}\label{GMEq2}
  \dot x = f(t,x,u),
\end{equation}
\begin{equation}\label{GMEq3}
  x(0)=x_0.
\end{equation}

{\it Подзадача улучшения\/} состоит в следующем: имеется допустимый, но не оптимальный процесс $v_s=(x_s(\cdot), u_s(\cdot))$. Требуется улучшить его, т.е. найти другой допустимый процесс $v_{s+1}=(x_{s+1}, u_{s+1})$, для которого $I(v_{s+1})<I(v_s)$. На каждой итерации численного метода итеративного типа решается подзадача улучшения текущего процесса. В процессе решения задачи строится последовательность $\{v_s\}$, причем $I(v_{s+1})<I(v_s)$ при всех $s$, и, если минимум $I_*$ в задаче (\ref{GMEq1})--(\ref{GMEq3}) конечен, то существует предел $I_-=\lim_{s\to\infty} I(v_s)$. Равенство $I_-=I_*$ в данной статье не исследуется. Для исследования этого вопроса можно применять, например, двойственный метод  \cite{BulatovKrotov2009}.

Глобальный метод улучшения управления основан на специальном преобразовании исходной задачи (\ref{GMEq1})--(\ref{GMEq3}) к эквивалентной ей задаче с другими $f^0$ и $F$. Это преобразование состоит в следующем: пусть функция $\varphi(t,x)$ непрерывно дифференцируема. Для любого процесса $v=(x(\cdot),u(\cdot))$, удовлетворяющего условиям (\ref{GMEq2}), (\ref{GMEq3}), в силу формулы Ньютона-Лейбница
\begin{equation}\label{GMEq4}
  \int\limits_0^T
  \left(
    \frac{\partial\varphi}{\partial t}(t,x(t))+
    \sum\limits_{i=1}^n
       \frac{\partial\varphi}{\partial x^i}(t,x(t))f^i(t,x(t),u(t))
  \right)\,dt -
  \varphi(T,x(T)) + \varphi(0,x_0) = 0.
\end{equation}
Целевой функционал для преобразованной задачи (\ref{GMEq1})--(\ref{GMEq3}) получается путем вычитания (\ref{GMEq4}) из (\ref{GMEq1}). Он имеет вид
\begin{displaymath}
L_\varphi (x;u) = -\int\limits_0^T R(t,x(t),u(t))\,dt + G(x(T)),
\end{displaymath}
где
\begin{equation}\label{GMDefR}
R(t,x,u)=-f^0(t,x,u)+ \frac{\partial\varphi}{\partial
t}(t,x) + \sum\limits_{i=1}^n \frac{\partial\varphi}{\partial
x^i}(t,x)f^i(t,x,u),
\end{equation}
\begin{equation}\label{GMDefG}
G(x)=F(x)+\varphi(T,x)-\varphi(0,x_0).
\end{equation}
Это преобразование выполняется на каждом шаге глобального метода для специально подобранной функции $\varphi$.

Опишем последовательность вычислительных операций, составляющих итерацию глобального метода. Начальное приближение $v_0=(x_0(\cdot), u_0(\cdot))$ задается произвольно. Пусть процесс $v_s=(x_s(\cdot), u_s(\cdot))$ уже рассчитан на $s$-й итерации.

Находим вспомогательную функцию $\varphi_s(t,x)$ такую, что
\begin{equation}\label{GMEq5}
  R_s(t,x_s(t),u_s(t))=\min\limits_{x} R_s(t,x,u_s(t)) \qquad (0<t<T),
\end{equation}
\begin{equation}\label{GMEq6}
  G_s(x_s(T))=\max\limits_x G_s(x).
\end{equation}
Необходимые условия (\ref{GMEq5}) и (\ref{GMEq6}):
\begin{equation}\label{GMEq5a}
  \frac{\partial}{\partial x} R_s(t,x_s(t),u_s(t))=0,
\end{equation}
\begin{equation}\label{GMEq6a}
  \frac{\partial}{\partial x} G_s(x_s(T))=0.
\end{equation}
Эти условия гарантируют выполнение (\ref{GMEq5}) и (\ref{GMEq6}) при дополнительных предположениях о выпуклости функции (\ref{GMDefR}) и вогнутости функции (\ref{GMDefG}).

Находим функцию
\begin{displaymath}
  \tilde u_s (t,x)\in\Argmax\limits_u R_s (t,x,u),
\end{displaymath}
решая задачу Коши
\begin{equation}\label{s3zk1}
  \dot x = f(t,x, \tilde u_s(t,x)),
\end{equation}
\begin{equation}\label{s3zk2}
  x(0) = x_0,
\end{equation}
находим траекторию $x_{s+1} (\cdot)$ и управление $u_{s+1} (t) = \tilde u_s (t,x_{s+1}(t))$.

Выражение
\begin{displaymath}
  I(v_s)-I(v_{s+1}) =
      \int\limits_0^T
      \left( R_s (t,x_{s+1},u_{s+1}) - R_s (t, x_{s+1}, u_s) \right)\,dt +
\end{displaymath}
\begin{displaymath}
      + \int\limits_0^T
      \left( R_s (t,x_{s+1},u_s) - R_s (t, x_s, u_s) \right)\,dt
      + G_s (x_s(T)) - G_s (x_{s+1}(T))
\end{displaymath}
положительно в силу (\ref{GMEq5}), (\ref{GMEq6}) и соотношения
\begin{equation}\label{umax}
  u_{s+1} (t) = \tilde u_s (t,x_{s+1}(t))\in\Argmax\limits_u R_s (t,x_{s+1}(t),u)
\end{equation}
при условии, что улучшаемый процесс $v_s$ не удовлетворяет принципу максимума. Таким образом, описанный алгоритм приводит к улучшению текущего процесса.

\section{Особенности реализации глобального метода
улучшения управления в исследуемой задаче}\label{HGM}

Рассмотрим особенности реализации глобального метода. Функцию $\varphi_s$ строим в виде $\varphi_s(t,x)=\psi_{s,1}(t)x^1 + \psi_{s,2}(t)x^2 +\ldots+\psi_{s,n}(t)x^{n} = (\psi_s(t),x)$. В этом случае
\begin{equation}\label{HGM_R}
  R_s(t,x,u) =
    \sum\limits_{i=1}^{n} \dot\psi_{s,i}(t)x^i +
    \sum\limits_{i=1}^{n} \psi_{s,i}(t)
           \sum\limits_{k=1}^{n} (a_{ik}+ub_{ik})x^k - \beta u^2;
\end{equation}
\begin{equation}\label{HGM_G}
  G_s(x) =
    \sum\limits_{i=1}^{n} \psi_{s,i}(T)x^i -
    \sum\limits_{i=1}^{n} \psi_{s,i}(0)x^i_0 +
           \sum\limits_{i,k=1}^{n} l_{ik}x^i x^k.
\end{equation}
Так как функция (\ref{HGM_R}) линейна по $x$, то из (\ref{GMEq5a}) следует (\ref{GMEq5}). Функция (\ref{HGM_G}) зависит от $x$ квадратично, поэтому (\ref{GMEq6}) следует из (\ref{GMEq6a}) при дополнительном предположении о вогнутости этой функции по пространственным переменным, т.е. о неположительной определенности её матрицы $L$. Так как все слагаемые в (\ref{HGM_G}), кроме последнего, линейны, то в этом случае функция (\ref{HGM_G}) вогнута, что обеспечивает возможность перехода от (\ref{GMEq6a}) к (\ref{GMEq6}). Таким образом, данный метод улучшения применим в общем случае (задача I) при дополнительном условии, что матрица $L$ неположительно определена. 

В задаче II мы можем добиться вогнутости функции (\ref{HGM_G}) путем замены исходной задачи на аналогичную задачу с матрицей $L_\alpha = L-\alpha I$, где $\alpha>0$ достаточно велико, или путём замены $L$ на ортопроектор на ортогональное дополнение (см. п.4-5 раздела \ref{SectPropsOfSystems}).

{\it Первый шаг} глобального метода - нахождение функций $\psi_i(\cdot)$ из условий (\ref{GMEq5a}), (\ref{GMEq6a}) - сводится к решению двойственной задачи Коши с терминальным условием.
\begin{displaymath}
  \frac{\partial R_s}{\partial x^m} = \dot\psi_{s,m} + \sum\limits_{i=1}^{n} (a_{im}+u_s b_{im})\psi_{s,i} = 0
  \qquad (\dot\psi_s=-(A^\transposed+u_s B^\transposed)\psi_s),
\end{displaymath}
\begin{displaymath}
  \frac{\partial G_s}{\partial x^m} = \psi_{s,m}(T)+2 \sum\limits_{i=1}^{n} l_{im}x^i  = 0
  \qquad (\psi_s(T)=-2Lx).
\end{displaymath}

Имеем согласно (\ref{HGM_R}):
\begin{displaymath}
  R_s(t,x,u) = \sum\limits_{i=1}^{n}\psi_{s,i}(t)x^i +
  \sum\limits_{i=1}^{n}
  \sum\limits_{k=1}^{n}
  \psi_{s,i}(t)a_{ik}x^k + u{\cal K} (\psi_s(t),x) - \beta u^2,
\end{displaymath}
где
\begin{displaymath}
  {\cal K} (\psi,x) =   \sum\limits_{i=1}^{n}
  \sum\limits_{k=1}^{n}
  \psi_i b_{ik} x^k = (Bx,\psi).
\end{displaymath}
Таким образом, при $\beta=0$
\begin{equation}\label{utilde}
  \tilde u_s(t, x) = \left\{
    \begin{array}{cc}
          \nu         &   \text{если ${\cal K} (\psi_s(t),x)>0$;} \\{}
           \text{любое из } [-\nu, \nu]     &    \text{если ${\cal K} (\psi_s(t),x)=0$;} \\{}
          -\nu          &    \text{если ${\cal K} (\psi_s(t),x)<0$,}
    \end{array}
    \right.
\end{equation}
а при $\beta>0$
\begin{displaymath}
  \tilde u_s(t, x) = \left\{
    \begin{array}{cc}
          -\nu         &   \text{если $\frac{1}{2\beta}{\cal K} (\psi_s(t),x) \leqslant -\nu$;} \\{}
           \frac{1}{2\beta}{\cal K} (\psi_s(t),x)     &    \text{если $-\nu<\frac{1}{2\beta}{\cal K} (\psi_s(t),x)<\nu$;} \\{}
          \nu          &    \text{если $\frac{1}{2\beta}{\cal K} (\psi_s(t),x) \geqslant \nu$.}
    \end{array}
    \right.
\end{displaymath}

{\it Второй шаг} сводится к решению задачи Коши (\ref{s3zk1}), (\ref{s3zk2}). Так как при $\beta>0$ правая часть системы (\ref{s3zk1}) непрерывна, то в этом случае решение задачи Коши теоретических проблем не составляет. Обратимся к случаю $\beta=0$. Основная трудность здесь состоит в том, что правая часть системы (\ref{s3zk1}) может иметь поверхности разрыва в фазовом пространстве. Задача (\ref{s3zk1}), (\ref{s3zk2}) решается поэтапно, каждый этап состоит в решении отдельной задачи Коши на участке непрерывности правой части системы (\ref{s3zk1}). Конечная точка траектории, найденная на очередном этапе, берется в качестве начального условия на следующем этапе. Обозначим через $\tau$ начальный момент времени для очередного этапа, а через $\xi_\tau = x_{s+1}(\tau)$ - соответствующую начальную точку траектории. Первоначально $\tau=0$, $\xi_\tau=x_0$. Обозначим $k(t)={\cal K} (\psi_s(t), x_{s+1}(t))$. Если $k(\tau)>0$, то решается задача Коши
\begin{equation}\label{ZKPlus}
   \dot x=f(t,x,\nu), \quad x(\tau)=\xi_\tau,
\end{equation}
а если $k(\tau)<0$, то решается задача
\begin{equation}\label{ZKMinus}
   \dot x=f(t,x,-\nu), \quad x(\tau)=\xi_\tau.
\end{equation}
Если $k(\tau)=0$, то выбор управления $u(t)$ при $t\geqslant\tau$ становится, вообще говоря, неоднозначным. Управление $u(t)\equiv \nu$ соответствует (\ref{utilde}), если на решении задачи  (\ref{ZKPlus}) функция $k(t)$ неотрицательна. Так как $k(\tau)=0$, то для неотрицательности $k(t)$ на некотором промежутке $(\tau, \tau+\varepsilon)$ достаточно, чтобы
\begin{equation}\label{dKdtPlus}
  \frac{d}{dt} k(\tau) =   (B(A+ \nu B) \xi_\tau, \psi_s(\tau)) - ((A+u_s (\tau) B)B \xi_\tau, \psi_s(\tau))>0.
\end{equation}
Управление $u(t)\equiv -\nu$ соответствует (\ref{utilde}), если на решении задачи  (\ref{ZKMinus}) функция $k(t)$ неположительна. Для этого достаточно, чтобы
\begin{equation}\label{dKdtMinus}
  \frac{d}{dt} k(\tau) =   (B(A- \nu B) \xi_\tau, \psi_s(\tau)) - ((A+u_s (\tau) B)B \xi_\tau, \psi_s(\tau))<0.
\end{equation}
Условия (\ref{dKdtPlus}) и (\ref{dKdtMinus}) не связаны между собой. Если выполняется только одно из этих условий, то выбор управления становится однозначным. Если выполняются оба условия, то произвольно выбирается одно из управлений $u(t)\equiv \nu$ или $u(t)\equiv -\nu$. Численное интегрирование соответствующей задачи Коши прекращается в тот момент, когда $k(t)=0$. Конечная точка найденной траектории выбирается в качестве $\xi_\tau$ на следующем этапе.

Если же не выполняется ни одно из условий (\ref{dKdtPlus}), (\ref{dKdtMinus}), то в качестве управления на текущем этапе выбирается управление $\hat u(\cdot)$, обеспечивающее $k(t)\equiv 0$ при $t\geqslant\tau$. По аналогии с принципом максимума назовем это управление {\it особым\/}. Из $dk/dt=0$ получаем
\begin{displaymath}
  u_{s+1}(t)=\hat u(t) = u_s(t) + \frac{((AB-BA) x_{s+1},\psi_s)}{(B^2 x_{s+1}, \psi_s)}.
\end{displaymath}
Таким образом, если на текущем этапе используется особое управление, то решается задача Коши
\begin{equation}\label{ZKDeg}
   \dot x=f\left(t,x,u_s(t) + \frac{((AB-BA)x,\psi_s)}{(B^2x, \psi_s)}\right), \quad x(\tau)=\xi_\tau.
\end{equation}
Решение задачи (\ref{ZKDeg}) прекращается в тот момент времени $t=\hat\tau>\tau$, когда будут выполнены условия, обеспечивающие возможность перехода с особого режима на управление $u(t)\equiv -\nu$ или $u(t)\equiv \nu$. Данные условия должны гарантировать выполнение (\ref{umax}) на некотором промежутке  $t\in(\hat\tau,\hat\tau+\varepsilon)$.

Рассмотрим один из возможных способов задания таких условий. Пусть $\hat\xi=x_{s+1}(\hat\tau)$ - конечная точка участка особого режима. Для возможности перехода с особого режима на управление $u(t)\equiv \nu$ необходимо, чтобы при всех $t$ из некоторого промежутка $t\in(\hat\tau,\hat\tau+\varepsilon)$ выполнялось неравенство ${\cal K} (\psi_s(t),x^+(t))\geqslant 0$, где $x^+(\cdot)$ - решение задачи
\begin{displaymath}
    \dot x^+=f(t,x^+,\nu), \quad x^+(\hat\tau)=\hat\xi.
\end{displaymath}
Так как ${\cal K} (\psi_s(\hat\tau),\hat\xi)= 0$, то для этого достаточно, чтобы $d{\cal K} (\psi_s(\hat\tau),\hat\xi)/dt> 0$. Последнее неравенство можно записать в виде
\begin{displaymath}
  ((BA-AB)\hat\xi, \psi_s(\hat\tau))+(\nu-u_s(\hat\tau))(B^2\hat\xi, \psi_s(\hat\tau))>0.
\end{displaymath}
Аналогичное условие, обеспечивающее возможность перехода с особого режима на управление $u(t)\equiv -\nu$, имеет вид
\begin{displaymath}
  ((BA-AB)\hat\xi, \psi_s(\hat\tau))-(\nu+u_s(\hat\tau))(B^2\hat\xi, \psi_s(\hat\tau))<0.
\end{displaymath}

Описанная процедура повторяется до тех пор, пока $\tau<T$.

Возможен другой способ реализации 3 шага, в котором не используется особый режим. Управление $u_{s+1}(\cdot)$ задается равенством
\begin{displaymath}
  \tilde u_{s+1}(t) = \left\{
    \begin{array}{cc}
          \nu         &   \text{если ${\cal K} (\psi_s(t),x_{s+1}(t)) \geqslant 0$;} \\{}
          -\nu          &    \text{если ${\cal K} (\psi_s(t),x_{s+1}(t))<0$.}
    \end{array}
    \right.
\end{displaymath}
Если решение задачи (\ref{s3zk1}), (\ref{s3zk2}) не содержит участков особого режима, то соответствующее управление является кусочно постоянной функцией, принимающей значения $-\nu$ и $\nu$, при этом количество переключений невелико и реализация численного решения задачи трудностей не представляет. При наличии особых режимов количество переключений управления становится очень большим, и с уменьшением шага интегрирования неограниченно возрастает. Так как управление с большим числом переключений малопригодно для практических целей, то участки с большим количеством переключений, на которых величина $k(t)$ близка к $0$, следует заменить на участки соответствующего особого режима, при этом управление (которое уже не будет константой) определяется из условия тождественного равенства нулю выражения ${\cal K} (\psi_s(t),x(t))$.

\section{Примеры}

В этом разделе описаны результаты численного решения задачи II различными итерационными методами в системе Matlab.

В силу наличия динамического инварианта абсолютный минимум квадратичной формы (\ref{MainProblem1}) на сфере $\|x\|=\|x_0\|$ является нижней оценкой для абсолютного минимума в задаче II. В первом примере (см. ниже) точка абсолютного минимума квадратичной формы (\ref{MainProblem1}) не принадлежит области достижимости задачи. Если параметры $T$ и $\nu$ достаточно велики, то допустимая траектория задачи II может оказаться в точке абсолютного минимума квадратичной формы (\ref{MainProblem1}) в момент времени $T$. Эта ситуация рассмотрена во втором примере.

Пример 1. Рассмотрим задачу: $n=4$, $T=\frac{1}{2}$, $\nu=3$,
\begin{displaymath}
  I(v) = -x_1^2\left(\frac{1}{2}\right)
  - 2x_2^2\left(\frac{1}{2}\right)
  - 3x_3^2\left(\frac{1}{2}\right)
  - x_4^2\left(\frac{1}{2}\right)
  \to min,
\end{displaymath}
\begin{displaymath}
 \dot{x} = (A+Bu)x, \qquad x(0) = \left(
   \begin{array}{c}
      -1\\
      1\\
     -1\\
     1
   \end{array}
\right), 
\end{displaymath}
где
\begin{displaymath}
A=\left(
\begin{array}{cccc}
0 & 0 & 1 & -2\\
0 & 0 & -2 & -1\\
-1 & 2 & 0 & 0\\
2 & 1 & 0 & 0
\end{array}
\right),\qquad
 B=\left(
\begin{array}{cccc}
0 & 0 & -1 & 1\\
0 & 0 & 1 & 2\\
1 & -1 & 0 & 0\\
 -1 & -2 & 0 & 0
\end{array}
\right).
\end{displaymath}
В ходе расчетов отрезок $[0, \frac{1}{2}]$  был аппроксимирован равномерной сеткой с шагом $0,0005$ ($1001$ узел), при этом значения управления, траектории и сопряженной траектории в узлах хранились в массивах. Их значения в промежуточных точках аппроксимировались линейными сплайнами. При решении задачи Коши использовался метод Рунге-Кутта четвертого порядка.

В качестве начального управления было взято постоянное управление $u_0(\cdot) \equiv 0.3.$ Задача решалась двумя методами - глобальным и градиентным \cite[Гл.6]{MDekker1996}.

На практике широко распространен технический прием, связанный с введением в целевой функционал задачи регуляризующего слагаемого, которое улучшает сходимость градиентного метода. Этот прием оказался полезным и при решении задачи глобальным методом. Рассмотрим вместо целевого функционала задачи следующий функционал
\begin{equation}\label{Iv}
  I(v) = -x_1^2\left(\frac{1}{2}\right) - 2x_2^2\left(\frac{1}{2}\right) -3x_3^2\left(\frac{1}{2}\right)-x_4^2\left(\frac{1}{2}\right)
  +\alpha\int\limits_0^{1/2} \left(u\left(t\right)\right)^2dt.
\end{equation}
Как показывает практика, при $\alpha>0$ градиентный метод в задаче (\ref{Iv}),  (\ref{MainProblem2})--(\ref{MainProblem4}) работает лучше, чем при $\alpha=0$ (задача I). Если $\alpha$ достаточно мало, то решение задачи (\ref{Iv}),  (\ref{MainProblem2})--(\ref{MainProblem4}) мало отличается от решения задачи I.

В таблице 1 представлены результаты решения задачи тремя методами.

Пример 2. Рассмотрим задачу:  $n=4$, $T=5$, $\nu=1$,
\begin{displaymath}
  I(v) = -x_1^2\left(5\right) - 2x_2^2\left(5\right) -3x_3^2\left(5\right)-x_4^2\left(5\right)\to  min,
\end{displaymath}
\begin{displaymath}
\dot{x} = (A+Bu)x, \qquad x(0) = \left(
\begin{array}{c}
-1\\
1\\
-1\\
1
\end{array}
\right),   
\end{displaymath}
где
\begin{displaymath}
  A=\left(
\begin{array}{cccc}
0 & 0 & 1 & -2\\
0 & 0 & -2 & -1\\
-1 & 2 & 0 & 0\\
2 & 1 & 0 & 0
\end{array}
\right),\qquad
 B=\left(
\begin{array}{cccc}
0 & 0 & -1 & 1\\
0 & 0 & 1 & 2\\
1 & -1 & 0 & 0\\
 -1 & -2 & 0 & 0
\end{array}
\right).
\end{displaymath}

Как и в предыдущем примере, отрезок $[0,5]$ был аппроксимирован
равномерной сеткой с шагом 0,0005. В качестве начального управления было взято кусочно постоянное управление $u_0(t) = 1$ при $t \in [0,1)$; $u_0(t) = 0$ при $t \in [1,5]$. Задача решалась градиентным методом, глобальным методом, а также глобальным методом с регуляризацией. В таблице 2 представлены результаты расчетов.

Полученные результаты позволяют сделать наблюдение, что на первых итерациях глобальный метод позволяет быстрее приблизиться к оптимальному управлению, чем градиентный. Добавление квадратичного интегрального слагаемого несколько ухудшает этот показатель, но снижает трудоемкость вычислительных операций.

\section{Заключение}

Исследована задача оптимального управления динамическими системами, описываемыми линейными дифференциальными уравнениями с управляемыми коэффициентами и квадратичным критерием оптимальности   характерная для многих  прикладных задач, играющих существенную роль в современной теории управления. Выделены некоторые специальные свойства таких систем, существенные для синтеза оптимального управления. Особый теоретический и прикладной интерес представляют гамильтоновы системы этого класса, для которых характерны компактная представимость в комплексном виде и наличие динамического инварианта, в силу которого все допустимые траектории расположены на сфере фазового пространства, фиксированной начальными условиями. Для понтрягинских экстремалей  рассмотренного класса задач характерно наличие особых режимов управления. Предпочтительный аппарат  решения данной задачи - итеративные методы оптимизации управления, состоящие в последовательном решении задачи улучшения, а из этих методов – глобальный метод улучшения. Последний детализирован и обоснован здесь применительно к данной задаче.

\appendix

\begin{proofoftheorem}{\ref{Th1}}
Задача Коши (\ref{MainProblem2})--(\ref{MainProblem3}) эквивалентна интегральному уравнению
\begin{displaymath}
  x(t) = x_0 + \int\limits_0^t (A+u(\tau)B)x(\tau) d\tau.
\end{displaymath}
Следовательно,
\begin{displaymath}
  \|x(t)\| \leqslant \|x_0\|  + \int\limits_0^t \|A+u(\tau)B\| \|x(\tau)\| d\tau
  \leqslant \|x_0\|  + \gamma \int\limits_0^t \|x(\tau)\| d\tau.
\end{displaymath}
Применяя неравенство Гронуолла-Беллмана \cite[стр.112]{GronwallBellman} получаем правое неравенство (\ref{DomainEstimate}). Для доказательства левого неравенства (\ref{DomainEstimate}) в задаче (\ref{MainProblem2})--(\ref{MainProblem3}) сделаем замену переменной $\theta=T-t$ и проведем аналогичные рассуждения.
\end{proofoftheorem}

\newpage
Таблица 1.

\begin{tabular}{|r|r|r|p{2.5cm}|}
\hline Номер итерации & Глобальный метод & Градиентный метод &
Глобальный метод с регуляризацией\\
\hline 0 &  -5,7363 & -5,7363 & -5,7363\\ \hline 
 1 & -6,3037 & -5,9802 & -6,1123\\ \hline
 2 & -6,6344 & -6,1917 & -6,2232\\ \hline
 3 & -8,0515 & -6,3612 & -6,8746\\ \hline
 4 & -9,1681 & -6,4917 & -8,3118\\ \hline
 5 & -9,4845 & -6,6005 & -9,148\\ \hline
 6 & -9,5743 & -6,7307 & -9,2799\\ \hline
 7 & -9,6592 & -6,9959 & -9,3817\\ \hline
 8 & -9,7570 & -7,6729 & -9,4776\\ \hline
 9 & -9,8323 & -8,8565 & -9,5896\\ \hline
 10 & -9,8908 & -9,6505 & -9,633\\ \hline
\end{tabular}

\newpage
Таблица 2.

\begin{tabular}{|r|c|c|p{2.5cm}|}
\hline Номер итерации & Глобальный метод & Градиентный метод &
Глобальный метод с регуляризацией\\
\hline 0 & -4,7545 & -4,7545 & -4,7545\\
\hline 1 & -7,2566 & -5,0698 & -6,1123\\ \hline
 2 &  -9,1684 & -5,4504 & -8,4596\\  \hline
 3 & -10,4411 & -6,0221 & -9,8425\\  \hline
 4 & -11,7116 & -6,8924 & -11,173\\  \hline
 5 & -11,9188 & -8,3881 & -11,455\\  \hline
 6 & -11,9335 & -10,187 & -11,516\\  \hline
 7 & -11,9784 & -11,643 & -11,525\\  \hline
 8 & -11,9900 & -11,826 & -11,534\\  \hline
 9 & -11,9931 & -11,907 & -11,505\\  \hline
 10 & -11,9939 & -11,942 & -11,529 \\ \hline
\end{tabular}

\newpage
\section*{Реферат}

Исследуются характерные свойства линейных дифференциальных уравнений с управляемыми коэффициентами, их динамические инварианты и ограничения области достижимости, существенные для проблем управления, свойства понтрягинских экстремалей данного класса задач, численные методы оптимизации управления.

\bigskip
\bigskip
\bigskip

02.30.Xx  - Calculus of variations

02.30.Yy	  - Control theory

02.60.Pn  - Numerical optimization

\newpage

Батурина Ольга Владимировна

аспирант Института проблем управления РАН

bulatov@ipu.ru

(495)334-91-59

\vspace{3 cm}

Булатов Александр Вячеславович

к.ф.-м.н.

старший научный сотрудник Института проблем управления РАН

bulatov@ipu.ru

(495)334-91-59

\vspace{3 cm}

Кротов Вадим Федорович

д.т.н., профессор

заведующий лабораторией Института проблем управления РАН

vfkrotov@ipu.ru

(495)334-91-59


\newpage

\begin{thebibliography}{99}

\bibitem{Krotov2009Ait}
  {\sl Кротов В.Ф.} Управление квантовыми системами и некоторые идеи оптимального управления. // АиТ. 2009. N. 3. C. 15–23.

\bibitem{Pontryagin}
  {\sl Понтрягин Л.С., Болтянский В.Г., Гамкрелидзе Р.В., Мищенко Е.Ф.} Математическая теория оптимальных процессов. М.: Наука, 1976.

\bibitem{GradMethodGeneral}
  {\sl Федоренко Р.П.} Приближенное решение задач оптимального управления. М.: Наука, 1978.

\bibitem{MDekker1996}
  {\sl Krotov V.F.} Global Methods in Optimal Control Theory. Marcel Dekker Inc.: NY, 1996.

\bibitem{GronwallBellman}
  {\sl Демидович Б.П.} Лекции по математической теории устойчивости. М.: Изд-во МГУ, 1998.

\bibitem{BulatovKrotov2009}
  {\sl Булатов А.В., Кротов В.Ф.} О численном решении линейно-квадратичной задачи оптимального управления двойственным методом. // АиT. 2009. N. 6. C. 3–14.

\end{thebibliography}
\end{document}